\documentclass[12pt]{article}
\usepackage{amsthm}
\usepackage{mathrsfs}
\usepackage{graphicx}
\usepackage{amsmath}
\usepackage{amsfonts}
\usepackage{amssymb}
\usepackage{url}
\usepackage{hyperref}
\usepackage{dsfont} 
\usepackage{verbatim}
\usepackage[mathlines]{lineno}
\usepackage{color}
\usepackage{enumerate}
\usepackage{enumitem}
\usepackage{graphicx}
\definecolor{red}{rgb}{1,0,0}

\definecolor{pink}{rgb}{.9,.2,.7}

\definecolor{org}{rgb}{1,.7,0}

\definecolor{purp}{rgb}{.5,0,.5}

\definecolor{br}{rgb}{.5,.7,1}

\newtheorem{thm}{Theorem}[section]

\newtheorem{prop}[thm]{Proposition}
\newtheorem{cor}[thm]{Corollary}
\newtheorem{rem}[thm]{Remark}
\newtheorem{lem}[thm]{Lemma}
\newtheorem{ex}[thm]{Example}
\newtheorem{defn}[thm]{Definition}

\newtheorem{obs}[thm]{Observation}

\def\noi{\noindent}
\def\mtx#1{\begin{bmatrix} #1 \end{bmatrix}}

\newcommand{\ba}{\begin{array}}
\newcommand{\ea}{\end{array}}
\newcommand{\x}{\times}
\newcommand{\bit}{\begin{itemize}}
\newcommand{\eit}{\end{itemize}}
\newcommand{\ben}{\begin{enumerate}}
\newcommand{\een}{\end{enumerate}}
\newcommand{\bea}{\begin{eqnarray*}}
\newcommand{\eea}{\end{eqnarray*}}
\newcommand{\beq}{\begin{equation}}
\newcommand{\eeq}{\end{equation}}
\newcommand{\bpf}{\begin{proof}}
\newcommand{\epf}{\end{proof}}

\newcommand{\diag}{\operatorname{diag}}

\DeclareMathSymbol{\mlq}{\mathord}{operators}{``}
\DeclareMathSymbol{\mrq}{\mathord}{operators}{`'}

\newcommand{\OL}{\overline}
\newcommand{\rank}{\operatorname{rank}}
\newcommand{\sepr}{\operatorname{sepr}}
\newcommand{\epr}{\operatorname{epr}}
\newcommand{\pr}{\operatorname{pr}}

\newcommand{\negat}{\operatorname{neg}}

\usepackage{authblk} 

\begin{document}



\title{{\bf The signed enhanced principal rank characteristic sequence}}





\author{Xavier Mart\'{i}nez-Rivera\thanks{
Corresponding author.
Email: xaviermr@iastate.edu}
}



\affil{\textit{Department of Mathematics, Iowa State University, Ames, IA, USA}}


\maketitle

\begin{abstract}
The signed enhanced principal rank characteristic sequence (sepr-sequence) of an $n \times n$ Hermitian matrix is the sequence $t_1t_2 \cdots t_n$, where
$t_k$ is either $\tt A^*$, $\tt A^+$, $\tt A^-$, $\tt N$, $\tt S^*$, $\tt S^+$, or $\tt S^-$, based on the following criteria:
$t_k = \tt A^*$ if $B$ has both a positive and a negative order-$k$ principal minor, and each order-$k$ principal minor is nonzero.
$t_k = \tt A^+$ (respectively, $t_k = \tt A^-$) if each order-$k$ principal minor is positive (respectively, negative).
$t_k = \tt N$ if each order-$k$ principal minor is zero.
$t_k = \tt S^*$ if $B$ has each a positive, a negative, and a zero order-$k$ principal minor.
$t_k = \tt S^+$ (respectively, $t_k = \tt S^-$) if $B$ has both a zero and a nonzero order-$k$ principal minor, and each nonzero order-$k$ principal minor is positive (respectively, negative).
Such sequences provide more information than the $({\tt A,N,S})$ epr-sequence in the literature,
where the $k$th term is either
$\tt A$, $\tt N$, or $\tt S$ based on whether
all, none, or some (but not all) of the
order-$k$ principal minors of the matrix are nonzero.
Various sepr-sequences are shown to be
unattainable by Hermitian matrices.
In particular,
by applying Muir's law of extensible minors,
it is shown that subsequences such as $\tt A^*N$ and $\tt NA^*$ are prohibited in the sepr-sequence of a Hermitian matrix.
For Hermitian matrices of orders $n=1,2,3$, all attainable sepr-sequences are classified.
For real symmetric matrices, a complete characterization of the attainable sepr-sequences whose underlying epr-sequence contains $\tt ANA$ as a non-terminal subsequence is established.
\end{abstract}

\noi{\bf Keywords.}
Signed enhanced principal rank characteristic sequence;
enhanced principal rank characteristic sequence;
minor; rank; Hermitian matrix.\\

\noi{\bf AMS Subject Classifications.}
15A15, 15A03, 15B57

\section{Introduction}
$\null$
\indent
The \textit{principal minor assignment problem}, introduced in \cite{PMAP introduced}, asks the following question: Can we find an $n \times n$ matrix with prescribed principal minors?
As a simplification of the principal minor assignment problem, Brualdi et al.\ \cite{ORIGINAL} associated a sequence with a symmetric matrix,
which they defined as follows:
Given an $n \times n$ symmetric matrix $B$ over a
field $F$,
\textit{the principal rank characteristic sequence}
(abbreviated pr-sequence) of $B$ is defined as
$\pr(B) = r_0]r_1 \cdots r_n$, where, for $k \geq 1$,
   \begin{equation*}
      r_k =
         \begin{cases}
             1 &\text{if $B$ has a nonzero principal minor of order $k$, and}\\
             0 &\text{otherwise,}
         \end{cases}
   \end{equation*}
while $r_0 = 1$ if and only if $B$ has a
$0$ diagonal entry \cite{ORIGINAL}.
(The \textit{order} of a minor is $k$ if it is
the determinant of a $k \times k$ submatrix.)
We note that the original definition of the pr-sequence was for real symmetric, complex symmetric and Hermitian matrices only;
but Barrett et al.\ \cite{FIELDS} later extended it to symmetric matrices over any field.
In the context of the principal minor assignment problem,
the pr-sequence is somewhat limited,
as it only records the presence or absence of a full-rank principal submatrix of each possible order;
thus, in order to provide more insight,
the pr-sequence was ``enhanced'' by Butler et al.\ \cite{EPR} with the introduction of another sequence:
Given an $n \times n$ symmetric matrix $B$ over a
field $F$,
the \textit{enhanced principal rank characteristic sequence} (abbreviated epr-sequence) of $B$ is defined as
$\epr(B) = \ell_1\ell_2 \cdots \ell_n$, where
   \begin{equation*}
      \ell_k =
         \begin{cases}
             \tt{A} &\text{if all the principal minors of order $k$ are nonzero;}\\
             \tt{S} &\text{if some but not all the principal minors of order $k$ are nonzero;}\\
             \tt{N} &\text{if none of the principal minors of order $k$ are nonzero, i.e., all are zero.}
         \end{cases}
   \end{equation*}
There has been substantial work done on pr- and epr-sequences
(see \cite{{ORIGINAL}, {FIELDS}, {EPR}, {skew}, {XMR-Classif}, {EPR-Hermitian}, {XMR-Char 2}}, for example).
Here, we introduce a sequence that extends the pr- and epr-sequence, which we think remains tractable, while providing further help for working on the principal minor assignment problem for Hermitian matrices:

\begin{defn}\rm{
Let $B$ be a complex Hermitian matrix with $\epr(B) = \ell_1\ell_2 \cdots \ell_n$.
The \textit{signed enhanced principal rank characteristic sequence} (abbreviated sepr-sequence) of $B$  is
the sequence $\sepr(B) = t_1t_2 \cdots t_n$, where
   \begin{equation*}
      t_k =
         \begin{cases}
             \tt{A^*} &\text{if $\ell_k = \tt{A}$
             and $B$ has both a positive and a negative order-$k$ principal minor;}\\
             \tt{A^+} &\text{if each order-$k$ principal minor of $B$ is positive;}\\
             \tt{A^-} &\text{if each order-$k$ principal minor of $B$ is negative;}\\
             \tt{N} &\text{if each order-$k$ principal minor of $B$ is zero;}\\
             \tt{S^*} &\text{if $\ell_k = \tt{S}$
             and $B$ has both a positive and a negative order-$k$ principal minor;}\\
             \tt{S^+} &\text{if $\ell_k = \tt{S}$ and each order-$k$ principal minor of $B$ is nonnegative;}\\
             \tt{S^-} &\text{if $\ell_k = \tt{S}$ and each  order-$k$ principal minor of $B$ is nonpositive.}\\
         \end{cases}
   \end{equation*}
}
\end{defn}

Further motivation for studying pr-, epr- and sepr-sequences are the instances where the principal minors of a matrix are of interest;
as stated in \cite{PMAP1},
these instances include the detection of $P$-matrices in the study of the complementarity problem,
Cartan matrices in Lie algebras,
univalent differentiable mappings,
self-validating algorithms,
interval matrix analysis,
counting of spanning trees of a graph using
the Laplacian, $D$-nilpotent automorphisms,
and in the solvability of the inverse multiplicative eigenvalue problem
(see \cite{PMAP1} and the references therein).

Section \ref{Sect:tools} is devoted to developing some of the tools used to establish results in subsequent sections.
In Section \ref{Sect:Hermitian}, various sepr-sequences are shown to be unattainable by Hermitian matrices.
Section \ref{Sect:Classif} is devoted to providing a classification of the sepr-sequences of orders $n=1,2,3$ that can be attained by an $n \times n$ Hermitian matrix.
Finally, Section \ref{Sect: real symm} focuses on the sepr-sequences of real symmetric matrices, where a complete characterization of the sepr-sequences whose underlying epr-sequence contains $\tt ANA$ as a non-terminal subsequence is established.

For the rest of the paper, all matrices are Hermitian.
For any sepr-sequence $\sigma$,
the epr-sequence resulting from
removing the superscripts of each term in $\sigma$ is
called the \emph{underlying} epr-sequence of $\sigma$.
A (pr-, epr- or sepr-) sequence is said to be
\textit{attainable} by a class of matrices provided that there exists a matrix $B$ in the class that attains it;
otherwise, we say that it is \textit{unattainable} (by the given class).
A subsequence that does not appear in an attainable
sequence is \textit{prohibited}.
A sequence is said to have \textit{order} $n$ if
it consists of $n$ terms.
Given a sequence
$t_{i_{1}} t_{i_{2}} \cdots t_{i_{k}}$, the notation
$\overline{t_{i_{1}} t_{i_{2}} \cdots t_{i_{k}}}$
indicates that the sequence may be repeated as many
times as desired (or it may be omitted entirely).
Let $B=[b_{ij}]$ and let
$\alpha, \beta \subseteq \{1, 2, \dots, n\}$;
then the submatrix lying in rows indexed by $\alpha$,
and columns indexed by $\beta$,
is denoted by $B[\alpha, \beta]$;
if $\alpha = \beta$, then the principal submatrix
$B[\alpha, \alpha]$ is abbreviated to $B[\alpha]$.
The matrices $O_n$, $I_n$ and $J_n$ are the matrices
of order $n$
denoting the zero matrix, the identity matrix and
the all-$1$s matrix, respectively.
The block diagonal matrix with the matrices $B$ and $C$ on the diagonal 
(i.e., the direct sum of $B$ and $C$) is denoted by $B \oplus C$.

\subsection{\textit{Results cited}}
$\null$
\indent
This section lists results that will be cited frequently, which will be referenced by the assigned nomenclature (if any).
Each instance of $\cdots$ below is permitted to be empty.

\begin{prop}\label{prop:SN...A...}
{\rm{\cite[Proposition 2.5]{EPR}}}
No Hermitian matrix can have the epr-sequence
${\tt SN}\cdots {\tt A}\cdots$.
\end{prop}

\begin{cor}\label{cor:noNSA}
{\rm{\cite[Corollary 2.7]{EPR}} ({\tt NSA} Theorem.)}
No Hermitian matrix can have ${\tt NSA}$ in its epr-sequence.
Further, no Hermitian matrix can have the epr-sequence  $\cdots {\tt ASN}\cdots {\tt A}\cdots$.
\end{cor}

For an $n \times n$ matrix $B$ with a nonsingular principal submatrix $B[\alpha]$, recall that the
Schur complement of $B[\alpha]$ in $B$ is the matrix
$B/B[\alpha] :=
B[\alpha^c] - B[\alpha^c,\alpha](B[\alpha])^{-1}B[\alpha,\alpha^c]$,
where $\alpha^c = \{1,2, \dots, n\} \setminus \alpha$.

\begin{thm}
\label{schur}
{\rm{\cite[Theorem 1.10]{EPR-Hermitian}}}
{\rm (Schur Complement Theorem.)}
Suppose $B$ is an $n \times n$
Hermitian matrix with $\rank B=r$.
Let $B[\alpha]$ be a nonsingular principal
submatrix of $B$ with $|\alpha| = k \leq r$,
and let $C = B/B[\alpha]$.
Then the following results hold.
\begin{enumerate}
\item [$(i)$]\label{p1SC} $C$ is an
$(n-k)\times (n-k)$ Hermitian matrix. 
\item [$(ii)$]\label{p2SC} Assuming the indexing of $C$ is inherited from $B$, any principal minor of $C$
is given by 
\[ \det C[\gamma] = \det B[\gamma \cup \alpha]/ \det B[\alpha].\]
\item [$(iii)$]\label{p3SC} $\rank C = r-k$.
\end{enumerate}
\end{thm}

\begin{cor}\label{schurAN}
{\rm \cite[Corollary 1.11]{EPR-Hermitian}}
{\rm (Schur Complement Corollary.)}
Let $B$ be a Hermitian matrix,
let $\epr(B)=\ell_1 \ell_2 \cdots \ell_n$, and
let $B[\alpha]$ be a nonsingular
principal submatrix of $B$,
with $|\alpha| = k \leq \rank B$.
Let $C = B/B[\alpha]$
and $\epr(C)=\ell'_{1} \ell'_2 \cdots \ell'_{n-k}$.
Then, for $j=1, \dots, n-k$,
$\ell'_j=\ell_{j+k}$  if
$\ell_{j+k} \in \{{\tt A,N}\}$.
\end{cor}

In the interest of brevity,
the notation $B_I$ for $\det(B[I])$
in \cite{ORIGINAL} and \cite{XMR-Char 2} is adopted here
(when $I = \emptyset$, $B_{\emptyset}$ is
defined to have the value 1).
Moreover, when $I = \{i_1, i_2, \dots, i_k\}$,
$B_I$ is written as $B_{i_1 i_2 \cdots i_k}$.

Given a matrix $B$, the determinant of the
$2 \times 2$ principal submatrix $B[\{i,j\}]$ can be stated as an homogenous polynomial identity as follows:
\[
B_{\emptyset}B_{ij} =
B_iB_j - \det(B[\{i\}|\{j\}])\det(B[\{j\}|\{i\}]).
\]
The identity in the next result is already known,
and can be obtained by applying
Muir's law of extensible minors \cite{Muir} to
the above identity
(for a more recent treatment of this law,
the reader is referred to \cite{Brualdi & Schneider}).

\begin{rem}\label{Muir Law}
\rm
Let $n \geq 2$,
let $B$ be an $n \times n$ Hermitian matrix,
let $i,j \in \{1,2, \dots, n\}$ be distinct, and
let $I \subseteq \{1,2, \dots, n\} \setminus \{i,j\}$.
Then
\[
B_{I}B_{I \cup \{i,j\}} =
B_{ I \cup \{i\} }B_{ I \cup \{j\} } -
|\det(B[ I \cup \{i\}| I \cup \{j\}])|^2.
\]
\end{rem}

Remark \ref{Muir Law} will be invoked as ``Muir's law of extensible minors.''
It should be noted that the identity in Remark \ref{Muir Law} can be derived from 
Sylvester's identity \cite{Brualdi & Schneider} in the case when $B_I \neq 0$.

\section{The signed enhanced principal rank characteristic sequence}\label{Sect:tools}
$\null$
\indent
We begin this section with simple observations,
and conclude with results that will serve as tools
in establishing the results of subsequent sections.

\begin{obs}
The sepr-sequence of a Hermitian
matrix  must end in $\tt A^+$, $\tt A^-$ or \ $\tt N$.
\end{obs}

Given an sepr-sequence $t_1t_2 \cdots t_n$,
the {\it negative} of this sequence, denoted
$\negat(t_1t_2 \cdots t_n)$,
is the sequence resulting from
replacing ``+'' superscripts with ``-'' superscripts
in $t_1t_2 \cdots t_n$, and viceversa.
For example, the negative of the sequence
${\tt NS^-S^*A^*A^+}$ is
${\tt NS^+S^*A^*A^-}$. Given a matrix $B$, the $i$th term in its sepr-sequence (respectively, epr-sequence) is
$[\sepr(B)]_i$ (respectively, $[\epr(B)]_i$).

\begin{obs}\label{odd terms obs}
Let $B$ be an $n \times n$ Hermitian matrix, and
let $i$ be an integer with $1 \leq i \leq n$.
\ben
\item If $i$ is even,
then $[\sepr(-B)]_{i} = [\sepr(B)]_{i}$.

\item If $i$ is odd, then
$[\sepr(-B)]_{i}$ = $\negat( [\sepr(B)]_{i})$.
\een
\end{obs}

The following is immediate from
\cite[Theorem 2.3]{EPR}.

\begin{thm}\label{NN result}{\rm ({\tt NN} Theorem.)}
Suppose $B$ is a Hermitian matrix,
$\sepr(B) = t_1t_2 \cdots t_n$, and
$t_k = t_{k+1} = \tt{N}$ for some $k$.
Then $t_j = \tt{N}$ for all $j \geq k$.
\end{thm}

\begin{thm}\label{Inverse Thm}
{\rm (Inverse Theorem.)}
Suppose $B$ is a nonsingular Hermitian matrix.
\ben
\item [$(i)$] If
$\sepr(B) = t_1t_2 \cdots t_{n-1}\tt{A^+}$, then
$\sepr(B^{-1}) = t_{n-1}t_{n-2} \cdots t_1 \tt{A^+}$.

\item [$(ii)$] If
$\sepr(B) = t_1t_2 \cdots t_{n-1}\tt{A^-}$, then
$\sepr(B^{-1})=
\negat(t_{n-1}t_{n-2} \cdots t_1)\tt{A^-}$.
\een
\end{thm}

\bpf
Let $\alpha \subseteq \{1,2, \dots, n\}$ be nonempty.
By Jacobi's determinantal identity,
$\det B^{-1}[\alpha] = \det B(\alpha) / \det B$.
The desired conclusions are now immediate.
\epf

The next lemma can be proven by replicating the proof of 
\cite[Theorem 2.6 (4)]{EPR} (the proof is virtually identical).

\begin{lem}\label{Dichot Lemma}
Let $k$ and $n$ be integers with $1 \leq k < n$.
Suppose that each $k$-element subset of
$\{1,2, \dots, n\}$ is associated with
exactly one of two given properties,
and that not every pair of $k$-element subsets is
associated with the same property.
Then there exist distinct integers
$i,j \in \{1,2, \dots, n\}$,
and a $(k-1)$-element subset
$I \subseteq \{1,2, \dots, n\} \setminus \{i,j\}$,
such that $I \cup \{i\}$ and $I \cup \{j\}$
are not associated with the same property.
\end{lem}


\begin{lem}\label{A*}
Let $B$ be an $n \times n$ Hermitian matrix
with $[\sepr(B)]_k = \tt A^*$.
Then there exists a
$(k-1)$-element subset $I \subseteq \{1,2, \dots, n\}$,
and $i,j \in \{1,2, \dots, n\} \setminus I$, such that
$B_{I \cup \{i\}} >0$ and $B_{I \cup \{j\}}<0$.
\end{lem}

\bpf
Since $\sepr(B)$ cannot end in $\tt A^*$, $k<n$.
Then, as every $k$-element subset of $\{1,2, \dots ,n\}$
is associated with either a positive or a negative
order-$k$ principal minor, but not both,
the conclusion follows from Lemma \ref{Dichot Lemma}. \epf

\begin{thm}\label{Inheritance}
{\rm (Inheritance Theorem.)}
Let $B$ be an $n \times n$ Hermitian matrix,
$m \leq n$, and $1\le i \le m$.
\ben
\item  If $[\sepr(B)]_i={\tt N}$, then  $[\sepr(C)]_i={\tt N}$ for all $m\times m$ principal submatrices $C$.
\item  If  $[\sepr(B)]_i={\tt A^+}$, then  $[\sepr(C)]_i={\tt A^+}$ for all $m\times m$ principal submatrices $C$.
\item  If  $[\sepr(B)]_i={\tt A^-}$, then  $[\sepr(C)]_i={\tt A^-}$ for all $m\times m$ principal submatrices $C$.
\item  If $[\sepr(B)]_m={\tt A^*}$, then there exist $m\times m$ principal submatrices $C_{A^+}$ and $C_{A^-}$ of $B$ such that $[\epr(C_{A^+})]_m = {\tt A^+}$ and $[\sepr(C_{A^-})]_m = {\tt A^-}$.
\item  If $[\sepr(B)]_m={\tt S^+}$, then there exist $m\times m$ principal submatrices $C_{A^+}$ and $C_N$ of $B$ such that $[\epr(C_{A^+})]_m = {\tt A^+}$ and $[\sepr(C_N)]_m = {\tt N}$.
\item  If $[\sepr(B)]_m={\tt S^-}$, then there exist $m\times m$ principal submatrices $C_{A^-}$ and $C_N$ of $B$ such that $[\sepr(C_{A^-})]_m = {\tt A^-}$ and $[\sepr(C_N)]_m = {\tt N}$.
\item  If $[\sepr(B)]_m={\tt S^*}$, then there exist $m\times m$ principal submatrices $C_{A^+}$, $C_{A^-}$ and $C_N$ of $B$ such that $[\sepr(C_{A^+})]_m = {\tt A^+}$, $[\sepr(C_{A^-})]_m = {\tt A^-}$ and $[\sepr(C_N)]_m = {\tt N}$.
\item  If $i < m$ and $[\sepr(B)]_i = {\tt A^*}$, then there exists an $m \times m$ principal submatrix $C_{A^*}$ such that $[\sepr(C_{A^*})]_i ={\tt A^*}$.
\item  If $i < m$ and $[\sepr(B)]_i = {\tt S^+}$, then there exists an $m \times m$ principal submatrix $C_{S^+}$ such that $[\sepr(C_{S^+})]_i ={\tt S^+}$.
\item  If $i < m$ and $[\sepr(B)]_i = {\tt S^-}$, then there exists an $m \times m$ principal submatrix $C_{S^-}$ such that $[\sepr(C_{S^-})]_i ={\tt S^-}$.
\item  If $i < m$ and $[\sepr(B)]_i = {\tt S^*}$, then there exists an $m \times m$ principal submatrix $C_{S}$ such that
    $[\sepr(C_S)]_i \in \{\tt S^*, S^+, S^-\}$.
\item  If $i < m$ and $[\sepr(B)]_i = {\tt S^*}$, then there exists an $m \times m$ principal submatrix $C_{+}$ such that
    $[\sepr(C_+)]_i \in \{\tt A^*, S^*, S^+\}$.

\item  If $i < m$ and $[\sepr(B)]_i = {\tt S^*}$, then there exists an $m \times m$ principal submatrix $C_{-}$ such that
    $[\sepr(C_-)]_i \in \{\tt A^*, S^*, S^-\}$.
\een
\end{thm}

\bpf
(1)--(3): Statements (1)--(3) simply follow by noting that a principal submatrix of a principal submatrix, is also principal submatrix.

(4)--(7): If $[\sepr(B)]_m = \tt A^*$, then
$B$ contains an $m \times m$ principal submatrix with
positive determinant, say, $C_{A^{+}}$, as well as one with negative determinant, say, $C_{A^{-}}$;
these two matrices each have the desired sepr-sequence,
which establishes Statement (4).
Statements (5)--(7) are established in the same manner as Statement (4).

(8): By Lemma \ref{A*}, there exists an
($i-1$)-element subset $I \subseteq \{1,2, \dots, m\}$,
and $p,q \in \{1,2, \dots, m\}\setminus I$,
such that
$B_{I \cup \{p\}} >0$ and $B_{I \cup \{q\}}<0$.
Then, by arbitrarily adding $m-i-1$ indices to
$I \cup \{p,q\}$, to obtain an $m$-element subset $\alpha$,
one obtains the principal submatrix $B[\alpha]$,
for which $[\sepr(B[\alpha])]_i = \tt A^*$.

(9)--(11): these three statements are immediate
from \cite[Theorem 2.6]{EPR}.

(12) and (13): By hypothesis, there exists two lists of indices, say,
$p_1,p_2,\dots,p_k$ and $q_1,q_2,\dots,q_k$, such that
$B_{p_1,p_2,\dots,p_i} >0$ and $B[q_1,q_2,\dots,q_i] <0$.
Without loss of generality, we may assume that
these lists are ordered so that any common indices occur in the same position in each list.
Consider the following lists of indices.
\[
\begin{array}{c}
p_1,p_2,p_3,\ldots,p_k;\\
q_1,p_2,p_3,\ldots,p_k;\\
q_1,q_2,p_3,\ldots,p_k;\\
q_1,q_2,q_3,\ldots,p_k;\\
\cdots\\
q_1,q_2,q_3,\ldots,q_k.
\end{array}
\]
As one moves down these lists, one must
eventually encounter two consecutive lists satisfying
one of the following:
(i) One list corresponds with a positive principal minor, and the other corresponds with a negative principal minor;
(ii) one list corresponds with a positive principal minor,
and the other corresponds with a zero principal minor.
If every pair of lists does not satisfy (i) or (ii),
then each list corresponds with a positive principal
minor, which is a contradiction, since the last list corresponds with a zero minor.
Hence, as two consecutive lists differ in
at most one position, the union of these two (distinct) lists generates an index set of cardinality $i+1$;
then, by arbitrarily adding $m-i-1$ indices to
this index set, to obtain an $m$-element subset $\alpha$,
one obtains the principal submatrix $B[\alpha]$,
for which
$[\sepr(B[\alpha])]_i \in \{\tt S^*, S^+ \}$ if the two
lists used to generate $\alpha$ satisfy (ii), while
$[\sepr(B[\alpha])]_i \in \{\tt A^*, S^* \}$ if the two
lists satisfy (i).
Hence, with $C_+ = B[\alpha]$,
$[\sepr(C_+)]_i \in \{\tt A^*, S^*, S^+\}$.

Statement (13) is established in the same manner as (12).
\epf

Given an $n \times n$ Hermitian matrix $B$
whose sepr-sequence contains
$\tt S^+$ (respectively, $\tt S^-$) in
position $i$, by the Inheritance Theorem,
for all $m$ with $i< m <n$, this matrix must
contain  at least one $m \times m$
principal submatrix whose sepr-sequence
inherits $\tt S^+$ (respectively, $\tt S^-$)
in position $i$.
However, the next example reveals that $\tt S^*$
is not necessarily inherited.

\begin{ex}{\rm
The (Hermitian) matrix
\[
B =
\mtx{
 -1 & 2 &   i &   4 & 0   \\
  2 & 0 &   6 &   1 & 8   \\
 -i & 6 &   1 &   i & 1+i \\
  4 & 1 &  -i &  -1 & 1+i \\
  0 & 8 & 1-i & 1-i & 0    }
\]
has sepr-sequence $\tt S^* S^- S^* A^+ A^+$.
It is easily verified that none of the
sepr-sequences of the five $4 \times 4$
principal submatrices of $B$ inherit the $\tt S^*$ appearing in the third position.
}
\end{ex}

With the next result, we add an additional tool to
our arsenal for studying epr- and sepr-sequences,
which is analogous to that of
the inheritance of an $\tt S^+$, $\tt S^-$ or $\tt A^*$ by
a principal submatrix (see the Inheritance Theorem).

%

\begin{prop}\label{schurA*}
Let $B$ be a Hermitian matrix with
$\sepr(B)= t_1 t_2 \cdots t_n$.
Suppose $t_p \in \{\tt A^*, A^+, A^- \}$
and     $t_q = \tt A^*$,
where $1 \leq p < q <n$.
Then there exists a (nonsingular)
$p \times p$ principal submatrix
$B[\alpha]$ such that the
$(n-p)\x (n-p)$ (Hermitian) matrix $C = B/B[\alpha]$
with $\sepr(C)= t'_{1} t'_2\cdots t'_{n-p}$
has $t'_{q-p}= t_{q} = \tt A^*$.
\end{prop}

\bpf
By Lemma \ref{Dichot Lemma},
there exist distinct integers $i,j \in \{1,2,\dots,n\}$,
and a ($q-1$)-element subset
$I \subseteq \{1,2, \dots, n\} \setminus \{i,j\}$,
such that
$\det B[I \cup \{i\}] > 0$ and
$\det B[I \cup \{j\}] < 0$.
Let $\alpha \subseteq I$ be a $p$-element subset.
By hypothesis, $B[\alpha]$ is nonsingular.
Let $C=B/B[\alpha]$,
$\sepr(C)= t'_{1} t'_2 \cdots t'_{n-p}$ and
$\beta = I \setminus \alpha$.
By the Schur Complement Theorem,
$\det (C[\beta \cup \{i\}])$ and
$\det (C[\beta \cup \{j\}])$ have opposite signs.
Then, as $|\beta \cup \{i\}| = |\beta \cup \{j\}| = q-p$,
$t'_{q-p} \in \{\tt A^*, S^*\}$.
But, by the Schur Complement Corollary,
we must have
$t'_{q-p} = {\tt A^*}$, as desired.
\epf

\section{Sepr-sequences of Hermitian matrices}\label{Sect:Hermitian}
$\null$
\indent
With our attention confined to Hermitian matrices,
in this section, we establish restrictions for
the attainability of sepr-sequences.

\begin{prop}\label{basic}{\rm (Basic Proposition.)}
No Hermitian matrix can have any of the
following sepr-sequences.
\ben
\item\label{basicA*A+} $\tt{A^*A^+ \cdots}$;
\item\label{basicA*S+} $\tt{A^*S^+ \cdots}$;
\item\label{basicA*N} $\tt{A^*N \cdots}$;
\item\label{basicS*A+} $\tt{S^*A^+ \cdots}$;
\item\label{basicS*S+} $\tt{S^*S^+ \cdots}$;
\item\label{basicS*N} $\tt{S^*N \cdots}$;
\item\label{basicS+A+} $\tt{S^+A^+ \cdots}$;
\item\label{basicS-A+} $\tt{S^-A^+ \cdots}$;
\item\label{basicNA*} $\tt{NA^* \cdots}$;
\item\label{basicNA+} $\tt{NA^+ \cdots}$;
\item\label{basicNS*} $\tt{NS^* \cdots}$;
\item\label{basicNS+} $\tt{NS^+ \cdots}$.
\een
\end{prop}

\bpf
To see that the sequences
\ref{basicA*A+}, \ref{basicA*S+}, \ref{basicA*N},
\ref{basicS*A+}, \ref{basicS*S+} and \ref{basicS*N}
are prohibited, note that a Hermitian matrix containing both a positive and a negative diagonal entry,
must contain a negative principal minor of order 2.

The sequences
\ref{basicS+A+} and \ref{basicS-A+} are prohibited
because a Hermitian matrix with both a zero and a nonzero
diagonal entry, must contain a nonpositive principal minor of order 2.

Finally, the fact that the sequences
\ref{basicNA*}, \ref{basicNA+},
\ref{basicNS*}, and \ref{basicNS+} are prohibited
follows from the fact that the principal minors of
order 2 of a Hermitian matrix with zero diagonal are nonpositive.
\epf

Although the following result is surely known, 
we note that it follows easily from Muir's law of extensible minors.

\begin{lem}\label{same sign lemma}
Let $B$ be a Hermitian matrix with
$\rank(B) = r$.
Then all the nonzero principal minors of $B$ of
order $r$ have the same sign.
\end{lem}


\begin{cor}\label{A*NN & S*NN}
Neither the sepr-sequences $\tt A^*NN$, nor $\tt S^*NN$, can occur as a subsequence of the sepr-sequence of a Hermitian matrix.
\end{cor}

\bpf
Let $B$ be a Hermitian matrix with $\sepr(B)$
containing $\tt A^*NN$ or $\tt S^*NN$,
where the $\tt A^*$ or $\tt S^*$ of this subsequence
occurs in position $k$.
Then, by the $\tt NN$ Theorem, $\rank(B) = k$.
Hence, by Lemma \ref{same sign lemma},
every nonzero principal minor of order $k$ has
the same sign, which is a contradiction.
\epf

In order to generalize one of the assertions of
Corollary \ref{A*NN & S*NN}, we will now apply
Muir's law of extensible minors.

\begin{thm}\label{A*N and NA*}
Neither the sepr-sequence $\tt{A^*N}$, nor $\tt{NA^*}$, can occur as a subsequence of the sepr-sequence of a Hermitian matrix.
\end{thm}

\bpf
Let $B$ be a Hermitian matrix with
$\sepr(B) = t_1t_2 \cdots t_n$.
Suppose to the contrary that $t_kt_{k+1} = \tt{A^*N}$
for some $k$.
By Lemma \ref{A*}, there exists a
($k-1$)-element subset $I \subseteq \{1,2, \dots, n\}$,
and $i,j \in \{1,2, \dots, n\} \setminus I$, such that
$B_{I \cup \{i\}} >0$ and $B_{I \cup \{j\}}<0$;
hence, $B_{I \cup \{i\}}B_{I \cup \{j\}} < 0$.
Now, since $I$ does not contain $i$ and $j$,
and because $B$ is Hermitian,
Muir's law of extensible minors implies that,
\[
B_{I}B_{I \cup \{i,j\}} =
B_{ I \cup \{i\} }B_{ I \cup \{j\} } -
|\det(B[ I \cup \{i\}| I \cup \{j\}])|^2.
\]
Then, as $B_{I \cup \{i\}}B_{I \cup \{j\}} < 0$,
$B_{I}B_{I \cup \{i,j\}} < 0$,
implying that $B_{I \cup \{i,j\}} \neq 0$,
a contradiction to $t_{k+1} = \tt{N}$.
It follows that $\tt{A^*N}$ is prohibited.

To establish the final assertion,
we again proceed by contradiction.
Suppose $t_kt_{k+1} = \tt{NA^*}$
for some $k$.
By the Basic Proposition, $k \geq 2$.
Since $t_{k+1} = \tt A^*$,
By Lemma \ref{A*}, there exists a
$k$-element subset $I \subseteq \{1,2, \dots, n\}$,
and $i,j \in \{1,2, \dots, n\} \setminus I$, such that
$B_{I \cup \{i\}} >0$ and $B_{I \cup \{j\}}<0$;
hence, $B_{I \cup \{i\}}B_{I \cup \{j\}} < 0$.
Once again, we use the identity
\[
B_{I}B_{I \cup \{i,j\}} =
B_{ I \cup \{i\} }B_{ I \cup \{j\} } -
|\det(B[ I \cup \{i\}| I \cup \{j\}])|^2.
\]
Since $t_k = \tt{N}$, we have $B_{I} = 0$,
implying that
\[0=
B_{ I \cup \{i\} }B_{ I \cup \{j\} } -
|\det(B[ I \cup \{i\}| I \cup \{j\}])|^2 <0,
\]
a contradiction.
\epf

For the rest of the paper, we invoke
Theorem \ref{A*N and NA*} by just stating that
$\tt A^*N$ or $\tt NA^*$ is prohibited.

\begin{thm}\label{AXA}
For any \ $\tt X$,
if any of the sepr-sequences
$\tt A^+XA^+$ or $\tt A^-XA^-$
occurs in the sepr-sequence of a Hermitian matrix,
then $\tt X \in \{A^+, A^-\}$.
\end{thm}

\bpf
Let $B$ be a Hermitian matrix with
$\sepr(B) = t_1t_2 \cdots t_n$.
Suppose
$t_k = t_{k+2} = \tt{A^+}$ or $t_k = t_{k+2} = \tt{A^-}$,
for some $k$ with $1 \leq k \leq n-2$.
Suppose to the contrary that
$t_{k+1} \neq \tt{A^+}$ and $t_{k+1} \neq \tt{A^-}$.
Let $\epr(B) = \ell_1\ell_2 \cdots \ell_n$.
Let $i,j \in \{1,2,\dots, n\}$, with $i \neq j$, and
$I \subseteq \{1,2, \dots, n\} \setminus \{i,j\}$,
where $|I| = k$.
Since $B$ is Hermitian,
Muir's law of extensible minors implies that
\[
B_{I}B_{I \cup \{i,j\}} =
B_{ I \cup \{i\} }B_{ I \cup \{j\} } -
|\det(B[ I \cup \{i\}| I \cup \{j\}])|^2.
\]
By hypothesis, $B_{I}B_{I \cup \{i,j\}} >0$.
Then, as $i$, $j$ and $I$ were arbitrary,
we must have
$B_{ I \cup \{i\} }B_{ I \cup \{j\} } >0$
whenever $I \subseteq \{1,2, \dots, n\} \setminus \{i,j\}$
and $i,j \in \{1,2,\dots, n\}$ are distinct
(otherwise, the expression on the right side of the above identity would be nonpositive).
It follows that $\ell_{k+1} = \tt A$.
By hypothesis, $t_{k+1} = \tt{A^*}$.
By Lemma \ref{A*}, there exists a
$k$-element subset $I \subseteq \{1,2, \dots, n\}$,
and $i,j \in \{1,2, \dots, n\} \setminus I$, such that
$B_{I \cup \{i\}} >0$ and $B_{I \cup \{j\}}<0$;
hence, $B_{I \cup \{i\}}B_{I \cup \{j\}} < 0$,
which is a contradiction to the above argument.
\epf

Theorem \ref{AXA} raises the following question:
Can the subsequences
$\tt{A^+A^+A^+}$, $\tt{A^+A^-A^+}$,
$\tt{A^-A^-A^-}$, $\tt{A^-A^+A^-}$ occur in the
sepr-sequence of a Hermitian matrix?
In Section \ref{Sect:Classif}, we demonstrate that
the answer is affirmative.

%

\begin{prop}\label{S+S+... and S-S+... implies singular}
For \ $\tt X \in \{ A^*, A^+, A^-\}$,
the sepr-sequences \
$\tt S^+S^+ \cdots X \cdots$ and \
$\tt S^-S^+ \cdots X \cdots$
are prohibited for Hermitian matrices.
\end{prop}

\bpf
Let $B = [b_{ij}]$ be an $n \times n$
Hermitian matrix with sepr-sequence
$\tt{S^+S^+} \cdots$ or $\tt{S^-S^+} \cdots$.
Without loss of generality,
we may assume that $b_{11} = 0$.
Let $j \in \{2,3, \dots, n\}$.
By hypothesis, the order-2 principal minor
$B_{1j} = b_{11}b_{jj} - |b_{1j}|^2 = -|b_{1j}|^2$
is nonnegative, implying that $b_{ij} = 0$.
Since $j$ was arbitrary, it follows that
the first row of $B$ is zero, implying that $B$ is
singular.
We conclude that a Hermitian matrix with sepr-sequence
$\tt{S^+S^+} \cdots$, or $\tt{S^-S^+} \cdots$,
is singular.

Now, suppose to the contrary that $B$ has sepr-sequence
$\tt S^+S^+ \cdots X \cdots$ or
$\tt S^-S^+ \cdots X \cdots$,
where $\tt X \in \{A^*, A^+, A^-\}$ occurs in
position $k$.
By the Inheritance Theorem,
$B$ has a nonsingular $k \times k$
principal submatrix with sepr-sequence
$\tt S^+Y \cdots$ or
$\tt S^-Y \cdots$,
where $\tt Y \in \{A^+, S^+, N\}$.
It follows from Proposition \ref{prop:SN...A...} and
the Basic Proposition that
$\tt Y = S^+$, a contradiction to
the assertion in the previous paragraph.
\epf

\begin{cor}\label{SSA} None of the following sepr-sequences can
occur as a subsequence of the sepr-sequence of
a Hermitian matrix.
\ben
\item  ${\tt S^+S^*A^+}$;
\item  ${\tt S^-S^*A^-}$;
\item  ${\tt S^+S^+A^+}$;
\item  ${\tt S^-S^-A^-}$;
\item  ${\tt S^+S^-A^+}$;
\item  ${\tt S^-S^+A^-}$.
\een
\end{cor}

\bpf
Let $B$ be a Hermitian matrix with
$\sepr(B) = t_1t_2 \cdots t_n$, where $n\geq 3$.
Let $k \in \{1,2, \dots, n-2 \}$.
We proceed by contradiction.

(1):
Suppose that
$t_k t_{k+1} t_{k+2} = {\tt S^+S^*A^+}$.
By the Inheritance Theorem,
$B$ contains a $(k+2) \x (k+2)$ principal submatrix $B'$
whose sepr-sequence ends with $\tt XYA^+$,
where $\tt X \in \{A^+, S^+, N \}$ and
      $\tt Y \in \{S^*, S^+, S^- \}$.
By the Inverse Theorem,
$\sepr((B')^{-1}) = \tt YX \cdots A^+$,
which contradicts
the Basic Proposition or
Proposition \ref{prop:SN...A...} or
Proposition \ref{S+S+... and S-S+... implies singular}.
It follows that ${\tt S^+S^*A^+}$ is prohibited.

(2):
Suppose that
$t_k t_{k+1} t_{k+2} = {\tt S^-S^*A^-}$.
By the Inheritance Theorem,
$B$ contains a $(k+2) \x (k+2)$ principal submatrix $B'$
whose sepr-sequence ends with $\tt XYA^-$,
where $\tt X \in \{A^-, S^-, N \}$ and
      $\tt Y \in \{S^*, S^+, S^- \}$.
By the Inverse Theorem,
$\sepr((B')^{-1}) = \tt \negat(YX) \cdots A^-$,
which contradicts
the Basic Proposition or
Proposition \ref{prop:SN...A...} or
Proposition \ref{S+S+... and S-S+... implies singular}.
Hence, ${\tt S^-S^*A^-}$ is prohibited.

(3):
Suppose that $t_k t_{k+1} t_{k+2} = {\tt S^+S^+A^+}$.
By the Inheritance Theorem,
$B$ contains a $(k+2)\times (k+2)$
principal submatrix whose sepr-sequence ends with
$\tt XS^+A^+$, where $\tt X \in \{A^+, S^+, N\}$.
By the Inverse Theorem,
$\sepr(B^{-1}) = \tt S^+ X \cdots A^+$.
Then, as $\tt X \in \{A^+, S^+, N\}$,
we have a contradiction to
the Basic Proposition or
Proposition \ref{prop:SN...A...} or
Proposition \ref{S+S+... and S-S+... implies singular}.
It follows that ${\tt S^+S^+A^+}$ is prohibited.

(4):
Suppose that $t_k t_{k+1} t_{k+2} = {\tt S^-S^-A^-}$.
By the Inheritance Theorem,
$B$ contains a $(k+2)\times (k+2)$
principal submatrix whose sepr-sequence ends with
$\tt XS^-A^-$, where $\tt X \in \{A^-, S^-, N\}$.
By the Inverse Theorem,
$\sepr(B^{-1})
= \tt \negat(S^- X) \cdots A^-
= \tt S^+\negat(X) \cdots A^-$,
which contradicts
the Basic Proposition or
Proposition \ref{prop:SN...A...} or
Proposition \ref{S+S+... and S-S+... implies singular}.

(5) and (6):
If any of
${\tt S^+S^-A^+}$ or ${\tt S^-S^+A^-}$
was a subsequence of $\sepr(B)$, then
applying Observation \ref{odd terms obs} to
$\sepr(B)$ would contradict items (3) or (4) above.
\epf

\begin{prop}\label{X+A*Y+ and X-A*Y-}
None of the following sepr-sequences can occur as a subsequence of the sepr-sequence of a Hermitian matrix.

\ben
\item $\tt A^+A^*S^+ $;
\item $\tt A^-A^*S^- $;
\item $\tt S^+A^*A^+ $;
\item $\tt S^-A^*A^- $;
\item $\tt S^+A^*S^+ $;
\item $\tt S^-A^*S^- $.
\een

\end{prop}

\bpf
Let $B$ be a Hermitian matrix containing one of the sequences (1)--(6) in positions $k-1, k, k+1$.
By Lemma \ref{A*}, there exists a
$(k-1)$-element subset $I \subseteq \{1,2, \dots, n\}$,
and $i,j \in \{1,2, \dots, n\} \setminus I$, such that
$B_{I \cup \{i\}} >0$ and $B_{I \cup \{j\}}<0$;
hence, $B_{I \cup \{i\}}B_{I \cup \{j\}}<0$.
But, since $B_{I}B_{I \cup \{i,j\}} \geq 0$ by hypothesis,
the identity
\[
B_{I}B_{I \cup \{i,j\}} =
B_{ I \cup \{i\} }B_{ I \cup \{j\} } -
|\det(B[ I \cup \{i\}| I \cup \{j\}])|^2
\]
leads to a contradiction.
\epf

\begin{prop}\label{SNA}
Let $B$ be a Hermitian matrix with
$\epr(B) = \ell_1\ell_2 \cdots \ell_n$ and
$\sepr(B) = t_1t_2 \cdots t_n$.
Suppose $\ell_k \ell_{k+1} \ell_{k+2} = \tt{SNA}$.
Then $t_kt_{k+1}t_{k+2} = \tt{S^+NA^-}$
or   $t_kt_{k+1}t_{k+2} = \tt{S^-NA^+}$.
\end{prop}

\bpf
Suppose to the contrary that
$t_kt_{k+1}t_{k+2} \neq \tt{S^+NA^-}$ or
$t_kt_{k+1}t_{k+2} \neq \tt{S^-NA^+}$.
Since $\tt{NA^*}$ is prohibited,
$t_kt_{k+1}t_{k+2} \in
\{\tt{S^*NA^+, S^+NA^+,
      S^*NA^-, S^-NA^-}\}$.
First, suppose
$t_kt_{k+1}t_{k+2} \in \{\tt{S^*NA^+, S^+NA^+}\}$.
By the Inheritance Theorem, $B$ has a
$(k+2) \times (k+2)$ principal submatrix $C_+$ with
$\sepr(C_+) \in
\{\tt \cdots A^*NA^+, \cdots S^*NA^+, \cdots S^+NA^+ \}$.
It now follows from the Inverse Theorem that
$\sepr(C_+^{-1}) \in
\{\tt NA^* \cdots, NS^* \cdots, NS^+  \cdots \}$,
which contradicts the Basic Proposition.

Finally, suppose
$t_kt_{k+1}t_{k+2} \in \{\tt{S^*NA^-, S^-NA^-}\}$.
By the Inheritance Theorem, $B$ has a
$(k+2) \times (k+2)$ principal submatrix $C_-$ with
$\sepr(C_-) \in
\{\tt \cdots A^*NA^-, \cdots S^*NA^-, \cdots S^-NA^- \}$.
Hence, by the Inverse Theorem,
$\sepr(C_-^{-1}) \in
\{\tt \negat(NA^*) \cdots,
      \negat(NS^*) \cdots,
      \negat(NS^-) \cdots \}$,
which contradicts the Basic Proposition.
\epf

A result analogous to Theorem \ref{AXA}
follows from Corollary \ref{SSA}, and
Propositions \ref{X+A*Y+ and X-A*Y-} and \ref{SNA}.

\begin{thm}\label{SXA}
For any $\tt X$,
if any of the sepr-sequences
$\tt S^+XA^+$  or $\tt S^-XA^-$
occurs in the sepr-sequence of a Hermitian matrix,
then $\tt X \in \{A^+, A^-\}$.
\end{thm}

\begin{prop}\label{...ASS...A...}
For $\tt X \in \{A^*, A^+, A^- \}$,
the following sepr-sequences are prohibited for
Hermitian matrices.

\ben
\item $\tt \cdots A^+S^*S^+ \cdots X \cdots$;
\item $\tt \cdots A^+S^+S^+ \cdots X \cdots$;
\item $\tt \cdots A^+S^-S^+ \cdots X \cdots$;
\item $\tt \cdots A^-S^*S^- \cdots X \cdots$;
\item $\tt \cdots A^-S^-S^- \cdots X \cdots$;
\item $\tt \cdots A^-S^+S^- \cdots X \cdots$.
\een
\end{prop}

\bpf
Let $\tt X \in \{ A^*, A^+, A^-\}$, and
let $B$ be a Hermitian matrix.
We first discard the sequences (1)--(3) simultaneously,
and then the sequences (4)--(6).

(1)--(3):
Suppose to the contrary that
$\sepr(B) = \tt \cdots A^+YS^+ \cdots X \cdots$,
where $\tt Y \in  \{ S^*, S^+, S^- \}$
and $\tt X$ occurs in position $k$.
By the Inheritance Theorem, $B$ has a
nonsingular $k \times k$ principal submatrix $B'$
whose sepr-sequence contains $\tt A^+WZ$,
where $\tt W \in  \{ S^*, S^+, S^- \}$ and
$\tt Z \in  \{ A^+, S^+, N \}$.
By Theorem \ref{AXA}, $\tt Z \neq A^+$.
Since $B'$ is nonsingular,
the $\tt NSA$ Theorem implies that $\tt Z \neq N$.
It follows that $\tt Z = S^+$.
Hence, $B'$ contains $\tt A^+WS^+$.
Then, as $B'$ is nonsingular,
the Inverse Theorem implies that
$\sepr((B')^{-1})$ contains one of
the prohibited sequences
$\tt S^+WA^+$ and $\tt \negat(S^+WA^+) = S^-\negat(W)A^-$,
which contradicts Corollary \ref{SSA}.

(4)--(6):
Suppose to the contrary that
$\sepr(B) = \tt \cdots A^-YS^- \cdots X \cdots$,
where $\tt Y \in  \{ S^*, S^+, S^- \}$
and $\tt X$ occurs in position of $k$.
By the Inheritance Theorem, $B$ has a
nonsingular $k \times k$ principal submatrix $B'$
whose sepr-sequence contains $\tt A^-WZ$,
where $\tt W \in  \{ S^*, S^+, S^-  \}$ and
$\tt Z \in  \{ A^-, S^-, N \}$.
By Theorem \ref{AXA}, $\tt Z \neq A^-$.
Since $B'$ is nonsingular,
the $\tt NSA$ Theorem implies that $\tt Z \neq N$.
It follows that $\tt Z = S^-$.
Hence, $B'$ contains $\tt A^-WS^-$.
Then, as $B'$ is nonsingular,
the Inverse Theorem implies that
$\sepr((B')^{-1})$ contains one of
the prohibited sequences
$\tt S^-WA^-$ and $\tt \negat(S^-WA^-) = S^+\negat(W)A^+$,
which again contradicts Corollary \ref{SSA}.
\epf

\begin{prop}\label{ANS}
Let $B$ be a Hermitian matrix with
$\epr(B) = \ell_1\ell_2 \cdots \ell_n$ and
$\sepr(B) = t_1t_2 \cdots t_n$.
Suppose $\ell_k \ell_{k+1} \ell_{k+2} = \tt{ANS}$.
Then $t_kt_{k+1}t_{k+2} = \tt{A^+NS^-}$
or   $t_kt_{k+1}t_{k+2} = \tt{A^-NS^+}$.
\end{prop}

\bpf
Suppose to the contrary that
$t_kt_{k+1}t_{k+2} \neq \tt{A^+NS^-}$ and   $t_kt_{k+1}t_{k+2} \neq \tt{A^-NS^+}$.
Clearly, $t_{k+1} = \tt{N}$.
Since $\tt A^*N$ is prohibited,
$t_k \neq \tt{A^*}$.
Hence, $t_kt_{k+1}t_{k+2} \in
\{ \tt{A^+NS^*, A^+NS^+,A^-NS^*,A^-NS^- }  \}$.
In each case,
by the Inheritance Theorem,
$B$ has a $(k+2) \times (k+2) $ principal submatrix $B'$
with $\sepr(B') = \tt{A^+NA^+}$
or   $\sepr(B') = \tt{A^-NA^-}$,
a contradiction to Theorem \ref{AXA}.
\epf

A natural question to answer is, does a result
analogous to Theorems \ref{AXA} and \ref{SXA} hold
for subsequences of the form
$\tt A^+XS^+$ and $\tt A^-XS^-$?
In other words, are the sequences
$\tt A^+XS^+$ and $\tt A^-XS^-$ prohibited in
the sepr-sequence of a Hermitian matrix when
$\tt X \notin \{A^+, A^-\}$?
The answer is negative:
An $n \times n$ positive semidefinite matrix with
nonzero diagonal and rank $n-1$,
and containing principal minors of order
$2$ and $3$ that are equal to zero,
serves as a counterexample.
A simple example is $B = I_3 \oplus J_2$.
Also, observe that $-B$ begins with $\tt A^-S^+S^-$
(see Observation \ref{odd terms obs}).
With that being said, a relatively similar result to
Theorems \ref{AXA} and \ref{SXA} can still be obtained,
which is an immediate consequence of Propositions
\ref{X+A*Y+ and X-A*Y-},
\ref{...ASS...A...} and
\ref{ANS}:

\begin{thm}\label{AXS...A...}
For any $\tt X$ and for $\tt Y \in \{A^*, A^+, A^-\}$,
if any of the sepr-sequences \\
$\tt \cdots A^+XS^+ \cdots Y \cdots$ or
$\tt \cdots A^-XS^- \cdots Y \cdots$
is attainable by a Hermitian matrix,
then  $\tt X\in \{A^+, A^-\}$.
\end{thm}

\begin{cor}
Any Hermitian matrix with an sepr-sequence containing
any of the following subsequences is singular.

\ben
\item $\tt A^+ S^* S^+$;
\item $\tt A^+ S^+ S^+$;
\item $\tt A^+ S^- S^+$;
\item $\tt A^- S^* S^-$;
\item $\tt A^- S^+ S^-$;
\item $\tt A^- S^- S^-$.
\een

\end{cor}

\begin{prop}\label{SNS}
Let $B$ be a Hermitian matrix with
$\epr(B) = \ell_1\ell_2 \cdots \ell_n$ and
$\sepr(B) = t_1t_2 \cdots t_n$.
Suppose $\ell_k \ell_{k+1} \ell_{k+2} = \tt{SNS}$
for some $k$ with $1 \leq k \leq n-2$.
Then $t_kt_{k+1}t_{k+2} = \tt{S^*NS^*}$,
or   $t_kt_{k+1}t_{k+2} = \tt{S^+NS^-}$,
or   $t_kt_{k+1}t_{k+2} = \tt{S^-NS^+}$.
\end{prop}

\bpf
We proceed by contradiction.
Suppose that  $t_kt_{k+1}t_{k+2} \notin
\{ \tt S^*NS^*, S^+NS^-, S^-NS^+ \}$.
Hence, $t_kt_{k+1}t_{k+2}$ is one of
the sequences in the set
\[
\{ \tt{S^*NS^+, S^*NS^-, S^+NS^*,
       S^+NS^+, S^-NS^*, S^-NS^-}\}.
\]
We examine four cases:

\textit{Case 1}: $t_kt_{k+1}t_{k+2} = \tt{S^*NS^+}$.
Let $B[\alpha]$ be a $k \times k$ nonsingular
principal submatrix with $\det (B[\alpha]) > 0$.
By the Schur Complement Theorem, $B/B[\alpha]$
is an $(n-k) \times (n-k)$ Hermitian matrix with
every diagonal entry equal to zero;
moreover, $\rank(B/B[\alpha]) =
\rank(B) - k \geq (k+2) - k = 2$,
implying that $B/B[\alpha]$ has a nonzero
principal minor of order 2,
say, $\det ((B/B[\alpha])[\{i,j\}])$.
Since $B/B[\alpha]$ has zero diagonal,
$\det ((B/B[\alpha])[\{i,j\}]) < 0$.
By the Schur Complement Theorem,
\[ \det ((B/B[\alpha])[\{i,j\}]) =
\det B[\{i,j\} \cup \alpha]/ \det B[\alpha].\]
Then, as $\det B[\alpha] > 0$,
$\det B[\{i,j\} \cup \alpha] < 0$,
a contradiction to $t_{k+2} = \tt S^+$.

\textit{Case 2}: $t_kt_{k+1}t_{k+2} = \tt{S^*NS^-}$.
Let $B[\alpha]$ be a $k \times k$ nonsingular
principal submatrix with $\det (B[\alpha]) < 0$.
Just as in Case 1, $B/B[\alpha]$ is an
$(n-k) \times (n-k)$ Hermitian matrix with zero diagonal,
and with a nonzero principal minor of order 2.
Let $\det ((B/B[\alpha])[\{i,j\}])$ be a nonzero
principal minor of order 2.
Since $B/B[\alpha]$ has zero diagonal,
$\det ((B/B[\alpha])[\{i,j\}]) < 0$.
As in Case 1, by the Schur Complement Theorem,
\[ \det ((B/B[\alpha])[\{i,j\}]) =
\det B[\{i,j\} \cup \alpha]/ \det B[\alpha].\]
Then, as $\det B[\alpha] < 0$,
$\det B[\{i,j\} \cup \alpha] > 0$,
a contradiction to $t_{k+2} = \tt S^-$.

\textit{Case 3}: $t_kt_{k+1}t_{k+2} \in
\{\tt{S^+NS^*, S^+NS^+}\}$.
By the Inheritance Theorem, $B$ has a
$(k+2) \times (k+2)$ principal submatrix $B'$
with $\sepr(B') = t'_1t'_2 \cdots t'_{k+2}$
having $t'_{k+1}t'_{k+2}  = \tt NA^+$ and
$t'_k \in \{\tt A^+, S^+, N\}$.
By the $\tt NN$ Theorem, $t'_k \neq \tt N$.
It follows that
$t'_k t'_{k+1}t'_{k+2} \in \{\tt A^+NA^+,S^+NA^+\}$,
which contradicts
Theorems \ref{AXA} and \ref{SXA}.

\textit{Case 4}: $t_kt_{k+1}t_{k+2} \in
\{\tt{S^-NS^*, S^-NS^-}\}$.
By the Inheritance Theorem, $B$ has a
$(k+2)\times (k+2)$ principal submatrix $B'$
with $\sepr(B') = t'_1t'_2 \cdots t'_{k+2}$ having
$t'_{k+1}t'_{k+2}  = \tt NA^-$ and
$t'_k \in \{\tt A^-, S^-, N\}$.
By the $\tt NN$ Theorem, $t'_k \neq \tt N$.
It follows that
$t'_k t'_{k+1}t'_{k+2} \in \{\tt A^-NA^-,S^-NA^-\}$,
a contradiction to
Theorems \ref{AXA} and \ref{SXA}.
\epf

\section{Sepr-sequences of order $n \leq 3$}\label{Sect:Classif}
$\null$
\indent
This section is devoted towards classifying all the sepr-sequences of orders $n = 1, 2, 3$ that can be
attained by Hermitian matrices.

For $n = 1$, it is obvious that the only attainable sepr-sequences are $\tt A^+, A^-$ and $\tt N$.

For $n=2$, there are a total of 21 sepr-sequences
ending in ${\tt A^+}$, ${\tt A^-}$ or ${\tt N}$;
of these, the 3 sequences that start with ${\tt S^*}$
are not attainable, since a matrix of order 2 only contains two diagonal entries.
Of the remaining 18 sequences,
${\tt A^*A^+}$,
${\tt A^*N}$,
${\tt S^+A^+}$,
${\tt S^-A^+}$ and
${\tt NA^+}$
are not attainable by the Basic Proposition.
That leaves 13 sequences, which constitute the
sepr-sequences of order $n=2$ that are attainable by
Hermitian matrices.
These 13 sequences are listed in
Table \ref{tab:order2}, where a Hermitian matrix attaining each sequence is provided;
in the case where the matrix provided is expressed as the negative of another matrix, its
sepr-sequence can be verified by applying
Observation \ref{odd terms obs} to the sepr-sequence of the corresponding matrix.

\begin{ex}\label{ex:n=2}
{\rm Matrices for Table \ref{tab:order2}:
}
\[
M_{{\tt A^*A^-}} =
\mtx{
1 &  1 \\
1 & -1 }\!\!, \
M_{{\tt S^+A^-}} =
\mtx{
1 & 1 \\
1 & 0 }.
\]
\end{ex}

\begin{table}[h!]\caption{All sepr-sequences of
order $n=2$  that are attainable by Hermitian matrices.   \label{tab:order2}}
\begin{center}{\scriptsize
 \begin{tabular}{|l|l|l|}\hline
\text{Sepr-sequence}\qquad$\null$  & \text{Hermitian  matrix}\qquad\qquad$\null$ & \text{Result}\qquad$\null$ \\
\hline
{$\tt A^*A^-$ }&$M_{\tt A^*A^-}$&Example \ref{ex:n=2}\\
\hline
{$\tt A^+A^+$ } & $I_2$    &     \\
\hline
{$\tt A^+A^-$ } & $2J_2 -I_2$  &   \\
\hline
{$\tt A^+N$ } &  $J_2$  &    \\
\hline
{$\tt A^-A^+$ } &  $-I_2$   &     \\
\hline
{$\tt A^-A^-$ } &  $-(2J_2 -I_2)$ &   \\
\hline
{$\tt A^-N$ } & $-J_2$  &    \\
\hline
{$\tt NA^-$ } & $J_2-I_2$  &   \\
\hline
{$\tt NN$ } &  $O_2$  &    \\
\hline
{$\tt S^+A^-$} & $M_{\tt S^+A^-}$ & Example \ref{ex:n=2}\\
\hline
{$\tt S^+N$} & $\diag(1,0)$ &     \\
\hline
{$\tt S^-A^-$}& $-M_{\tt S^+A^-}$&Example \ref{ex:n=2}\\
\hline
{$\tt S^-N$}& $\diag(-1,0)$ &     \\
\hline
\end{tabular}}
\end{center}
\end{table}

As just shown, the results developed before
this section sufficed to
decide the attainability of all the sepr-sequences of order $n=2$. However, for $n=3$,
there remain sequences unaccounted for.

\begin{prop}{\rm (Order-3 Proposition)}
For any $\tt X$,
the sepr-sequences \
$\tt S^* S^* X$,
$\tt S^* A^* X$,
$\tt A^* S^* X$ and \
$\tt X S^* N$
are prohibited as the sepr-sequence of a \
$3 \times 3$ Hermitian matrix.
\end{prop}

\bpf
To see why
$\tt S^* S^* X$ and $\tt S^* A^* X$ are prohibited,
observe that any $3 \times 3$ Hermitian matrix
whose sepr-sequence starts with $\tt S^*$ cannot contain
a positive principal minor of order 2,
since such a matrix does not contain two nonzero diagonal entries having the same sign.

To discard $\tt A^* S^* X$,
note that a $3 \times 3$ Hermitian matrix with
an sepr-sequence starting with $\tt A^*$ must
contain at least two
negative principal minors of order 2.
Then, as a $3 \times 3$ matrix contains only 3
principal minors of order 2,
a $3 \times 3$ Hermitian matrix cannot
have an sepr-sequence
starting with $\tt A^*S^*$,
since it cannot have both a zero and a positive
principal minor of order 2.

Finally, the fact that $\tt X S^* N$ is prohibited
follows from Lemma \ref{same sign lemma},
since a Hermitian matrix attaining this sequence
would have rank 2.
\epf

Since the underlying epr-sequence of an
attainable sepr-sequence must also be attainable,
to decide the attainability of
the sepr-sequences of order $3$,
we will take advantage of what is known about
the epr-sequences of $3 \times 3$
Hermitian matrices.
We proceed by first determining the sepr-sequences that are attainable by Hermitian matrices but not by real symmetric matrices, and then
we determine the remaining ones,
namely those that can be attained by real symmetric matrices.

It was established in \cite{EPR-Hermitian} that
$\tt NAN$ is the only epr-sequence of order $3$ that is attainable by a Hermitian matrix but not by a real symmetric matrix.
Since $\tt NA^*N$ and $\tt NA^+N$ are not attainable by a Hermitian matrix (because of the Basic Proposition),
the only sepr-sequence of order 3 that is attainable by a Hermitian matrix but not by a real symmetric matrix is $\tt NA^-N$.

It now remains to determine the sepr-sequences that are attainable by real symmetric matrices.
The epr-sequences of order $3$ that
are attainable by real symmetric matrices are listed in \cite[Table 3]{EPR}, which are
$\tt AAA$,
$\tt AAN$,
$\tt ANA$,
$\tt ANN$,
$\tt ASA$,
$\tt ASN$,
$\tt NAA$,
$\tt NNN$,
$\tt NSN$,
$\tt SAA$,
$\tt SAN$,
$\tt SNN$,
$\tt SSA$ and
$\tt SSN$.
Then, as an attainable sepr-sequence must end in
$\tt A^+$, $\tt A^-$ or $\tt N$,
by counting the sepr-sequences whose underlying
epr-sequence is one of those just listed,
we find that only $130$ sepr-sequences are
potentially attainable
(note that we are \textit{not}
counting the sequence $\tt NA^-N$ among
these 130 sequences,
since we are now only counting those that are
attainable by real symmetric matrices).
We now discard certain sequences from these 130 sequences,
and show that the remaining ones are all attainable.
The 3 sepr-sequences starting with $\tt A^*A^+$ are not attainable by the Basic Proposition;
that leaves 127 sequences.
The 10 sequences having one of the forms
$\tt A^+XA^+$ or $\tt A^-XA^-$,
with $\tt X \notin \{A^+, A^-\}$,
are not attainable by Theorem \ref{AXA};
that leaves 117 sequences.
The 11 sequences containing the prohibited subsequences $\tt A^*N$ and $\tt NA^*$ are discarded;
that leaves 106 sequences.
Of the remaining sequences
(which do not include the
already-discarded sequence $\tt S^*A^*N$),
13 are discarded by the Order-3 Proposition;
that leaves 93 sequences.
The 3 sequences of the form $\tt A^* S^+X$,
as well as
$\tt NA^+A^+$, $\tt NA^+A^-$ and $\tt NS^+N$,
are discarded by the Basic Proposition;
that leaves sequences  87 sequences.
The 9 sequences starting with $\tt XA^+$,
where $\tt X \in \{S^*, S^+, S^-\}$,
are discarded by the Basic Proposition;
that leaves 78 sequences.
The 8 sequences of the form
$\tt S^+XA^+$ and $\tt S^-XA^-$,
with $\tt X \in \{A^*, S^*, S^+, S^- \}$,
are discarded by
Theorem \ref{SXA};
that leaves 70 sequences.
The sequence $\tt S^*NN$ is discarded by the
Basic Proposition;
that leaves 69 sequences.
The sequences
$\tt S^+S^+A^-$ and $\tt S^-S^+A^+$
are discarded by
Proposition \ref{S+S+... and S-S+... implies singular};
that leaves 67 sequences.
Finally, the 3 sequences of the form $\tt S^*S^+X$ are discarded by the Basic Proposition;
that leaves 64 sequences,
which constitute the sepr-sequences of order $n=3$ that are attainable by real symmetric matrices.
By adding the sequence $\tt NA^-N$ to these 64 sequences, we obtain all the sepr-sequences that are attainable by Hermitian matrices;
these 65 sequences are listed in
Table \ref{tab:order3}, where a Hermitian matrix
attaining each sequence is provided.

\begin{ex}\label{ex:n=3}
{\rm
Matrices for Table \ref{tab:order3}:
{\small
\[
M_{\tt A^*A^-A^+} =
\mtx{
 1 & 2 & 2 \\
 2 & 1 & 2 \\
 2 & 2 & -1}\!\!, \
M_{\tt A^+A^+A^-} =
\mtx{
 1 & 1 & -1 \\
 1 & 2 & 1 \\
 -1 & 1 & 2}\!\!, \
M_{\tt A^+A^-A^+} =
\mtx{
 1 & 2 & 2 \\
 2 & 1 & 2 \\
 2 & 2 & 1}\!\!, \
\]
\[
M_{\tt A^+A^-A^-} =
\mtx{
 1 & 2 & -2 \\
 2 & 1 & 2 \\
 -2 & 2 & 1}\!\!, \
M_{\tt A^*A^-N} =
\mtx{
 1 & 2 & 0 \\
 2 & 1 & \sqrt{3} \\
 0 & \sqrt{3} & -1}\!\!, \
M_{\tt A^+A^+N} =
\mtx{
 2 & 1 & 1 \\
 1 & 2 & -1 \\
 1 & -1 & 2}\!\!, \
\]
\[
M_{\tt A^+A^-N} =
\mtx{
 1 & 2 & 2 \\
 2 & 1 & 7 \\
 2 & 7 & 1}\!\!, \
M_{\tt A^*S^-A^+} =
\mtx{
 -1 & -1 & 0 \\
 -1 & -1 & -1 \\
 0 & -1 & 1}\!\!, \
M_{\tt A^+S^*A^-} =
\mtx{
 1 & -2 & -4 \\
 -2 & 4 & 2 \\
 -4 & 2 & 4}\!\!, \
\]
\[
M_{\tt A^+S^+A^-} =
\mtx{
 1 & 1 & 0 \\
 1 & 1 & 1 \\
 0 & 1 & 1}\!\!, \
M_{\tt A^+S^-A^-} =
\mtx{
 1 & 1 & 2 \\
 1 & 1 & 3 \\
 2 & 3 & 1}\!\!, \
M_{\tt A^*S^-N} =
\mtx{
 -1 & 0 & 0 \\
 0 & 1 & 1 \\
 0 & 1 & 1}\!\!, \
\]
\[
M_{\tt A^+S^-N} =
\mtx{
 1 & 2 & 2 \\
 2 & 1 & 1 \\
 2 & 1 & 1}\!\!, \
M_{\tt NA^-N} =
\mtx{
  0 & i & 1 \\
 -i & 0 & 1 \\
  1 & 1 & 0}\!\!, \
M_{\tt NS^-N} =
\mtx{
 0 & 1 & 0 \\
 1 & 0 & 0 \\
 0 & 0 & 0}\!\!, \
\]
\[
M_{\tt S^+A^*A^-} =
\mtx{
 2 & 1 & 1 \\
 1 & 2 & 2 \\
 1 & 2 & 0}\!\!, \
M_{\tt S^+A^-A^+} =
\mtx{
 1 & 1 & 1 \\
 1 & 0 & 1 \\
 1 & 1 & 0}\!\!,\
M_{\tt S^+A^-A^-} =
\mtx{
 2 & 1 & 1 \\
 1 & 0 & 2 \\
 1 & 2 & 0}\!\!, \
\]
\[
M_{\tt S^*A^-N} =
\mtx{
 1 & 0 & 1 \\
 0 & -1 & 1 \\
 1 & 1 & 0}\!\!, \
M_{\tt S^+A^-N} =
\mtx{
 2 & 2 & 1 \\
 2 & 0 & 2 \\
 1 & 2 & 0}\!\!, \
M_{\tt S^*S^-A^+} =
\mtx{
 1 & 0 & 1 \\
 0 & -1 & 0 \\
 1 & 0 & 0}\!\!, \
\]
\[
M_{\tt S^+S^-A^-} =
\mtx{
 1 & 0 & 0 \\
 0 & 0 & 1 \\
 0 & 1 & 0}\!\!, \
M_{\tt S^+S^+N} =
\mtx{
 2 & 1 & 0 \\
 1 & 2 & 0 \\
 0 & 0 & 0}\!\!, \
M_{\tt S^+S^-N} =
\mtx{
 1 & 1 & 1 \\
 1 & 0 & 0 \\
 1 & 0 & 0}\!\!.
\]
}
}
\end{ex}

\begin{table}[h!]\caption{All sepr-sequences of
order $n=3$
that are attainable by
Hermitian matrices.}\label{tab:order3}
\begin{center}{\scriptsize
 \begin{tabular}{|l|l|l|}\hline
\text{Sepr-sequence}\qquad$\null$  & \text{Hermitian matrix}\qquad\qquad$\null$ & \text{Result}\qquad$\null$ \\ \hline
$\tt A^*A^*A^+$ & $\diag(1,-1,-1)$ & \quad$\null$\\

\hline
$\tt A^*A^*A^-$ & $\diag(-1,1,1)$  &  \\

\hline
$\tt A^*A^-A^+$ & $M_{\tt A^*A^-A^+}$    &
Example \ref{ex:n=3}   \\

\hline
$\tt A^*A^-A^-$ & $-M_{\tt A^*A^-A^+}$    &
Example \ref{ex:n=3}  \\

\hline
$\tt A^+A^*A^-$ & $(-M_{\tt A^*A^-A^+})^{-1}$ &
Example \ref{ex:n=3}  \\

\hline
$\tt A^+A^+A^+$ & $I_3$   &   \\

\hline
$\tt A^+A^+A^-$ & $M_{\tt A^+A^+A^-}$    &
Example \ref{ex:n=3}   \\

\hline
$\tt A^+A^-A^+$ & $M_{\tt A^+A^-A^+}$    &
Example \ref{ex:n=3}   \\

\hline
$\tt A^+A^-A^-$ & $M_{\tt A^+A^-A^-}$    &
Example \ref{ex:n=3}   \\

\hline
$\tt A^-A^*A^+$ & $-(-M_{\tt A^*A^-A^+})^{-1}$    &
Example \ref{ex:n=3}  \\

\hline
$\tt A^-A^+A^+$ & $-M_{\tt A^+A^+A^-}$    &
Example \ref{ex:n=3}  \\

\hline
$\tt A^-A^+A^-$ & $-I_3$    &   \\

\hline
$\tt A^-A^-A^+$ & $-M_{\tt A^+A^-A^-}$    &
Example \ref{ex:n=3}  \\

\hline
$\tt A^-A^-A^-$ & $-M_{\tt A^+A^-A^+}$    &
Example \ref{ex:n=3}  \\

\hline
$\tt A^*A^-N$ & $M_{\tt A^*A^-N}$ &
Example \ref{ex:n=3}   \\

\hline
$\tt A^+A^+N$ & $M_{\tt A^+A^+N}$   &
Example \ref{ex:n=3}  \\

\hline
$\tt A^+A^-N$ & $M_{\tt A^+A^-N}$ &
Example \ref{ex:n=3}  \\

\hline
$\tt A^-A^+N$ &$-M_{\tt A^+A^+N}$ &
Example \ref{ex:n=3}  \\

\hline
$\tt A^-A^-N$ & $-M_{\tt A^+A^-N}$ &
Example \ref{ex:n=3}  \\


\hline
$\tt A^+NA^-$ & $-(J_3 - 2I_3)$   &   \\

\hline
$\tt A^-NA^+$ & $J_3 - 2I_3$    &   \\


\hline
$\tt A^+NN$ & $J_3$    &   \\

\hline
$\tt A^-NN$ & $-J_3$   &   \\


\hline
$\tt A^*S^-A^+$& $M_{\tt A^*S^-A^+}$   &
Example \ref{ex:n=3}   \\

\hline
$\tt A^*S^-A^-$& $-M_{\tt A^*S^-A^+}$   &
Example \ref{ex:n=3}  \\

\hline
$\tt A^+S^*A^-$ & $M_{\tt A^+S^*A^-}$  &
Example \ref{ex:n=3}   \\

\hline
$\tt A^+S^+A^-$& $M_{\tt A^+S^+A^-}$   &
Example \ref{ex:n=3}   \\

\hline
$\tt A^+S^-A^-$ & $M_{\tt A^+S^-A^-}$   &
Example \ref{ex:n=3}   \\

\hline
$\tt A^-S^*A^+$ & $-M_{\tt A^+S^*A^-}$   &
Example \ref{ex:n=3}  \\

\hline
$\tt A^-S^+A^+$ & $-M_{\tt A^+S^+A^-}$   &
Example \ref{ex:n=3}  \\

\hline
$\tt A^-S^-A^+$ &  $-M_{\tt A^+S^-A^-}$  &
Example \ref{ex:n=3}  \\


\hline
$\tt A^*S^-N$ & $M_{\tt A^*S^-N}$   &
Example \ref{ex:n=3}   \\

\hline
$\tt A^+S^+N$ & $J_1 \oplus J_2$   &   \\

\hline
$\tt A^+S^-N$ & $M_{\tt A^+S^-N}$   &
Example \ref{ex:n=3}   \\

\hline
$\tt A^-S^+N$ & $-(J_1 \oplus J_2)$   &   \\

\hline
$\tt A^-S^-N$& $-M_{\tt A^+S^-N}$   &
Example \ref{ex:n=3}  \\


\hline
$\tt NA^-A^+$ & $J_3 - I_3$   &   \\

\hline
$\tt NA^-A^-$ & $-(J_3 - I_3)$   &   \\

\hline
$\tt NA^-N$ & $M_{\tt NA^-N}$    &
Example \ref{ex:n=3}   \\


\hline
$\tt NNN$     & $O_3$    &   \\

\hline
$\tt NS^-N$ & $M_{\tt NS^-N}$    &
Example \ref{ex:n=3}   \\


\hline
$\tt S^*A^-A^+$ & $(-M_{\tt A^+S^*A^-})^{-1}$  &
Example \ref{ex:n=3}  \\

\hline
$\tt S^*A^-A^-$    & $-(-M_{\tt A^+S^*A^-})^{-1}$    &
Example \ref{ex:n=3}  \\

\hline
$\tt S^+A^*A^-$   &  $M_{\tt S^+A^*A^-}$  &
Example \ref{ex:n=3}   \\

\hline
$\tt S^+A^-A^+$     &  $M_{\tt S^+A^-A^+}$   &
Example \ref{ex:n=3}   \\

\hline
$\tt S^+A^-A^-$    &  $M_{\tt S^+A^-A^-}$   &
Example \ref{ex:n=3}   \\

\hline
$\tt S^-A^*A^+$    & $-M_{\tt S^+A^*A^-}$    &
Example \ref{ex:n=3}  \\

\hline
$\tt S^-A^-A^+$    & $-M_{\tt S^+A^-A^-}$   &
Example \ref{ex:n=3}  \\

\hline
$\tt S^-A^-A^-$    & $-M_{\tt S^+A^-A^+}$   &
Example \ref{ex:n=3}  \\


\hline
$\tt S^*A^-N$ & $M_{\tt S^*A^-N}$    &
Example \ref{ex:n=3}   \\

\hline
$\tt S^+A^-N$ & $M_{\tt S^+A^-N}$    &
Example \ref{ex:n=3}   \\

\hline
$\tt S^-A^-N$ & $-M_{\tt S^+A^-N}$    &
Example \ref{ex:n=3}  \\


\hline
$\tt S^+NN$ & $J_1 \oplus O_2$   &   \\

\hline
$\tt S^-NN$ & $-(J_1 \oplus O_2)$   &   \\

\hline
$\tt S^*S^-A^+$ & $M_{\tt S^*S^-A^+}$    &
Example \ref{ex:n=3}   \\

\hline
$\tt S^*S^-A^-$ & $-M_{\tt S^*S^-A^+}$    &
Example \ref{ex:n=3}  \\

\hline
\end{tabular}}
\end{center}
\end{table}

\pagebreak

\noindent Table \ref{tab:order3} (continued): All sepr-sequences of 
order $n=3$  that are attainable by
Hermitian matrices.
\begin{center}{\scriptsize
 \begin{tabular}{|l|l|l|}\hline
\text{Sepr-sequence}\qquad$\null$  & \text{Hermitian matrix}\qquad\qquad$\null$ & \text{Result}\qquad$\null$ \\

\hline
$\tt S^+S^*A^-$ & $(-M_{\tt S^*S^-A^+})^{-1}$    &
Example \ref{ex:n=3}  \\


\hline
$\tt S^+S^-A^-$    & $M_{\tt S^+S^-A^-}$   &
Example \ref{ex:n=3}   \\

\hline
$\tt S^-S^*A^+$     & $-(-M_{\tt S^*S^-A^+})^{-1}$    &
Example \ref{ex:n=3}  \\

\hline
$\tt S^-S^-A^+$    & $-M_{\tt S^+S^-A^-}$   &
Example \ref{ex:n=3}  \\


\hline
$\tt S^*S^-N$    & $\diag(1,-1,0)$    &   \\

\hline
$\tt S^+S^+N$     & $M_{\tt S^+S^+N}$    &
Example \ref{ex:n=3}   \\

\hline
$\tt S^+S^-N$     & $M_{\tt S^+S^-N}$    &
Example \ref{ex:n=3}   \\

\hline
$\tt S^-S^+N$    & $-M_{\tt S^+S^+N}$    &
Example \ref{ex:n=3}  \\

\hline
$\tt S^-S^-N$    & $-M_{\tt S^+S^-N}$  &
Example \ref{ex:n=3}  \\

\hline
\end{tabular}}
\end{center}

\section{Sepr-sequences of real symmetric matrices}\label{Sect: real symm}
$\null$
\indent
This section focuses on real symmetric matrices,
and its main result is a complete characterization of the sepr-sequences whose underlying epr-sequence contains
$\tt ANA$ as a non-terminal subsequence (see Theorem \ref{ANA nonterminal characterization}).

\begin{prop}
For any $\tt X$, the sepr-sequence $\tt NXS^*N$
cannot occur in the sepr-sequence of a real symmetric matrix.
\end{prop}

\bpf
Let $B$ be a real symmetric matrix
with $\sepr(B)$ containing $\tt NXS^*N$,
where the penultimate term of this subsequence
occurs in position $k$.
By {\cite[Proposition 2.4]{XMR-Classif}},
$\rank(B) = k$.
It follows from Lemma \ref{same sign lemma}
that the nonzero principal minors of order $k$ of $B$
have the same sign, which contradicts our hypothesis.
\epf

%
%

\begin{prop}\label{ANA at start effect}
Let $B$ be a real symmetric matrix with
$\epr(B) = \ell_1 \ell_2 \cdots \ell_n$ and
$\sepr(B) = t_1 t_2 \cdots t_n$.
Suppose $\ell_1 \ell_2 \ell_3 = {\tt ANA}$.
Then
$t_1 t_2 t_3 = {\tt A^+NA^-}$ or
$t_1 t_2 t_3 = {\tt A^-NA^+}$.
Furthermore, the following hold.
\ben
\item If $t_{1} t_{2}t_{3} = \tt A^+NA^-$, then
$t_i = {\tt A^-}$ for $i \geq 4$.

\item If $t_{1} t_{2}t_{3} = \tt A^-NA^+$, then
$t_i = {\tt A^+}$ for odd $i \geq 4$, and
$t_j = {\tt A^-}$ for even $j \geq 4$.

\een
\end{prop}

\bpf
By \cite[Proposition 2.5]{XMR-Classif},
$B$ is conjugate by a nonsingular diagonal matrix to one of the matrices $\pm$($J_n - 2I_n$).
Since $\sepr(B)$ remains invariant under this type of conjugation, we may assume that $B = \pm$($J_n - 2I_n$).
It is now easy to check that
$t_1 t_2 t_3 = {\tt A^+NA^-}$ or
$t_1 t_2 t_3 = {\tt A^-NA^+}$.
We examine each case separately.

\textit{Case 1}: $t_1 t_2 t_3 = {\tt A^+NA^-}$.
Hence, $B =-$($J_n - 2I_n$).
Let $k$ be an integer with $4 \leq  k \leq n$.
Observe that any order-$k$ principal submatrix is of the
form $-$($J_k - 2I_k$);
hence, each order-$k$ principal minor is
$2^{k-1}(2-k) < 0$
(the eigenvalues of $-$($J_k - 2I_k$) are $2-k$ and $2$, with multiplicity $1$ and $k-1$, respectively).
It follows that $t_k = \tt A^-$.

\textit{Case 2}: $t_1 t_2 t_3 = {\tt A^-NA^+}$.
Hence, $B =J_n - 2I_n$.
The desired conclusion now follows
by applying Observation \ref{odd terms obs}
to the matrix $-B$, which, by Case 1, has
$\sepr(B) = \tt A^+NA^- \OL{A^-}$.
\epf

\begin{cor}\label{A-NA+A+ is prohibited}
The sepr-sequence $\tt A^-NA^+A^+$ cannot occur as a subsequence of the sepr-sequence of a real symmetric matrix.
\end{cor}

\bpf
If a real symmetric matrix existed with an sepr-sequence containing $\tt A^-NA^+A^+$, then, by the Inheritance Theorem, it would contain a principal submatrix whose sepr-sequence ends with $\tt A^-NA^+A^+$, and whose
inverse has sepr-sequence $\tt A^+NA^- \cdots A^+$
(see the Inverse Theorem), which would contradict
Proposition \ref{ANA at start effect}.
\epf

\begin{cor}\label{odd-even cor: A+NA-A+ & A-NA+A-}
Let $B$ be a real symmetric matrix with
$\sepr(B) = t_1 t_2 \cdots t_n$, and
let $k$ be an integer with $k \leq n-2$.
Then the following hold.
\ben
\item
If $t_{k-1} t_k t_{k+1} t_{k+2} = \tt A^+NA^-A^+$,
then $k$ is odd.

\item
If $t_{k-1} t_k t_{k+1} t_{k+2} = \tt A^-NA^+A^-$,
then $k$ is even.
\een
\end{cor}

\bpf
Suppose
$t_{k-1} t_k t_{k+1} t_{k+2} = \tt A^+NA^-A^+$.
If $k$ were even, then,
by Observation \ref{odd terms obs},
$\sepr(-B)$ would contain $\tt A^-NA^+A^+$,
which would contradict
Proposition \ref{A-NA+A+ is prohibited}.
That establishes Statement (1).
Statement (2) is proven similarly.
\epf

\begin{lem}\label{ANA lemma}
Let $B$ be a real symmetric matrix with
$\epr(B) = \ell_1 \ell_2 \cdots \ell_n$ and
$\sepr(B) = t_1 t_2 \cdots t_n$.
Suppose $\ell_{k-1} \ell_{k}\ell_{k+1} = \tt ANA$,
where $k \leq n-2$.
Then
$t_i \in \{ \tt A^+, A^- \}$ for all $i \neq k$
and
$t_{k+1} = \negat(t_{k-1})$.
\end{lem}

\bpf
If $k=2$, then all the conclusions are immediate from
Proposition \ref{ANA at start effect};
thus, we assume that $k \geq 3$.
By  \cite[Theorem 2.6]{XMR-Classif},
$t_i \in \{ \tt A^*, A^+, A^- \}$
for all $i \neq k$.
By Theorem \ref{AXA}, and because
$\tt A^*N$ and $\tt NA^*$ are prohibited,
$t_{k-1} t_k t_{k+1} = \tt A^+NA^-$ or
$t_{k-1} t_k t_{k+1} = \tt A^-NA^+$;
hence, $t_{k+1} = \negat(t_{k-1})$.
We now show by contradiction that
$t_i \neq \tt A^*$ for all $i \neq k$;
thus, suppose $t_j = \tt A^*$ for some $j \neq k$.
We proceed by examining two cases.

\textit{Case 1}: $j < k$.
Since $t_{k+2} \in \{\tt A^*, A^+,A^-\}$,
the Inheritance Theorem implies that
$B$ has a (necessarily nonsingular)
$(k+2) \times (k+2)$ principal submatrix $B'$ with
$\sepr(B') = \tt \cdots A^* \cdots XN\negat(X)Y$,
where $\tt X, Y \in \{A^+, A^- \}$.
By the Inverse Theorem, \\
$\sepr((B')^{-1}) = \tt Z N \negat(Z) \cdots A^* \cdots$,
where $\tt Z \in \{ \tt A^+, A^- \}$;
now observe that this contradicts
Proposition \ref{ANA at start effect}.

\textit{Case 2}: $j>k$.
Since $t_{k-2} \in \{\tt A^*, A^+,A^-\}$,
Proposition \ref{schurA*} implies that
there exists a (necessarily nonsingular)
$(k-2) \times (k-2)$ principal submatrix
$B[\alpha]$ such that the sepr-sequence of
$C = B/B[\alpha]$ has $\tt A^*$ in the
$(j-(k-2))$-th position.
By the Schur Complement Corollary,
$\epr(C) = \tt ANA \cdots$;
hence, by Proposition \ref{ANA at start effect},
$\sepr(C)$ does not contain $\tt A^*$,
which leads to a contradiction.
\epf

We are now in position to completely characterize all
the sepr-sequences that are attainable by real symmetric matrices and whose underlying epr-sequence contains
$\tt ANA$ as a non-terminal subsequence.

\begin{thm}\label{ANA nonterminal characterization}
Let $\sigma = t_1 t_2 \cdots t_n$ be
an sepr-sequence whose underlying epr-sequence is
$\ell_1 \ell_2 \cdots \ell_n$.
Suppose $\ell_{k-1} \ell_{k}\ell_{k+1} = \tt ANA$,
where $2 \leq k \leq n-2$.
Let $\alpha = \{1, \dots, n-1\} \setminus \{k-1, k\}$.
Then $\sigma$ is attainable by a real symmetric matrix
if and only if
one of the following holds.

\ben

\item
$\sigma = \tt \OL{A^+} A^+ N A^- A^- \OL{A^-}$;

\item $k$ is odd,
$t_{k-1} t_k t_{k+1} t_{k+2} = \tt A^+NA^-A^+$ and
$t_{i+1} = \negat(t_i)$ for all
$i \in \alpha$;

\item $k$ is even,
$t_{k-1} t_k t_{k+1} t_{k+2} = \tt A^-NA^+A^-$ and
$t_{i+1} = \negat(t_i)$ for all
$i \in \alpha$.

\een
\end{thm}

\bpf
First, we show that if any of Statements (1)--(3)
holds, then $\sigma$ is attainable.
Suppose (1) holds.
Let $B = -(J_n - kI_n)$.
We claim that $\sepr(B) = \sigma$.
Obviously, $[\sepr(B)]_1 = {\tt A^+} = t_1$.
Since every principal submatrix of
order $q \geq 2$ is of the form $-(J_q - kI_q)$,
each principal minor of order $q$ is
$k^{q-1}(k-q)$.
Hence,
$[\sepr(B)]_q = {\tt A^+} = t_q$ for $2 \leq q \leq k-1$,
$[\sepr(B)]_k = {\tt N} = t_k$, and
$[\sepr(B)]_q = {\tt A^-} = t_q$ for $k+1 \leq q \leq n$.
It follows that $\sepr(B) = \sigma$.
To show that $\sigma$ is attainable if
Statements (2) or (3) hold, let $C = J_n - kI_n$.
Note that
each principal minor of order $q \geq 2$ is
$(-k)^{q-1}(q-k) = (-k)^q(k-q)$.
It is now easy to check that $\sepr(C) = \sigma$,
with $\sigma$ as in Statement (2) or (3),
depending on the parity of $k$.

For the other direction,
suppose $\sigma$ is attainable by a
real symmetric matrix, say, $B$,
so that $\sepr(B) = \sigma = t_1t_2 \cdots t_n$.
By Lemma \ref{ANA lemma},
$t_i \in \{ \tt A^+, A^- \}$ for all $i \neq k$
and
$t_{k+1} = \negat(t_{k-1})$.
It follows that
$t_{k-1} t_k t_{k+1} = \tt A^+NA^-$ or
$t_{k-1} t_k t_{k+1} = \tt A^-NA^+$.
Since $\tt A^-NA^+A^+$ is prohibited by
Corollary \ref{A-NA+A+ is prohibited},
$t_{k-1} t_k t_{k+1} t_{k+2}$ must be either
$\tt A^+NA^-A^-$,
$\tt A^+NA^-A^+$ or
$\tt A^-NA^+A^-$.
We now examine these three possibilities in two cases.

\textit{Case i}: $t_{k-1} t_k t_{k+1} t_{k+2} = \tt A^+NA^-A^-$.
We now show that
$\sepr(B) = \tt \OL{A^+} A^+ N A^- A^- \OL{A^-}$.
We start by showing that $t_i = \tt A^+$ for
all $i \leq k-2$.
Suppose to the contrary that
there exists $j \leq k-2$ such that
$t_j = \tt A^-$.
By the Inheritance and Inverse Theorems,
the sepr-sequence of the inverse of
any (necessarily nonsingular) $(k+2) \times (k+2)$
principal submatrix of $B$ has the form
$\tt A^+ N A^- \cdots A^+ \cdots$,
which contradicts
Proposition \ref{ANA at start effect}.
We conclude that $t_i = \tt A^+$ for
all $i \leq k-2$.
Now we show that $t_i = \tt A^-$ for
all $i \geq k+3$.
Suppose to the contrary that there exists
$j \geq k+3$  such that $t_j = \tt A^+$.
Then,
as every principal minor of order $k-2$ is positive,
the Schur Complement Theorem and
the Schur Complement Corollary imply that
for any (necessarily nonsingular)
$(k-2) \times (k-2)$ principal submatrix $B[\alpha]$,
$\sepr(B/B[\alpha]) = \tt A^+ N A^- \cdots A^+ \cdots$, which contradicts
Proposition \ref{ANA at start effect}.
We conclude that
$\sepr(B) = \tt \OL{A^+} A^+ N A^- A^- \OL{A^-}$.
Then, as $\sigma = \sepr(B)$, Statement (1) holds.
Note that we have shown that if the sepr-sequence of
a real symmetric matrix contains $\tt A^+NA^-A^-$,
then its sepr-sequence must be
$\tt \OL{A^+} A^+ N A^- A^- \OL{A^-}$.

\textit{Case ii}:
$t_{k-1} t_k t_{k+1} t_{k+2} = \tt A^+NA^-A^+$ or
$t_{k-1} t_k t_{k+1} t_{k+2} = \tt A^-NA^+A^-$.
By Corollary \ref{odd-even cor: A+NA-A+ & A-NA+A-},
$k$ is odd if
$t_{k-1} t_k t_{k+1} t_{k+2} = \tt A^+NA^-A^+$,
and
$k$ is even if
$t_{k-1} t_k t_{k+1} t_{k+2} = \tt A^-NA^+A^-$.
Thus, it remains to show that
$t_{i+1} = \negat(t_i)$ for all
$i \in \alpha$, from which it would follow that
either Statement (2) or Statement (3) holds.
Suppose to the contrary that
$t_{i+1} \neq \negat(t_i)$ for some
$i \in \alpha$;
hence, $t_{i+1} = t_i$.
Let $\sepr(-B) = t'_1t'_2 \cdots t'_n$.
It follows from Observation \ref{odd terms obs} that
$t'_{k-1} t'_k t'_{k+1} t'_{k+2} = \tt A^+NA^-A^-$
and that
$t'_{i+1} = \negat(t'_i)$.
Now, observe that the last sentence at the end of Case i
implies that
$\sepr(-B) = \tt \OL{A^+} A^+ N A^- A^- \OL{A^-}$.
Hence, $t'_{i+1} = t'_i$.
Then, as $t'_{i+1} = \negat(t'_i)$, we must have
$t'_{i+1} = t'_i = \tt N$, a contradiction.
We conclude that it must be the case that
$t_{i+1} = \negat(t_i)$ for all
$i \in \alpha$, implying that either
Statement (2) or Statement (3) holds.
\epf

To see that Theorem \ref{ANA nonterminal characterization} cannot be generalized to Hermitian matrices,
the reader is referred to
\cite[Theorem 3.3]{EPR-Hermitian}.
Moreover, Theorem \ref{ANA nonterminal characterization} cannot be generalized to include the case when $\tt ANA$ occurs as a terminal sequence, since the epr-sequence $\tt SAANA$ is attainable by a real symmetric matrix (see \cite[Table 5]{EPR}).


\section*{Acknowledgements}
$\null$
\indent
The author expresses his gratitude to Dr. Leslie Hogben, 
for introducing him to the topic of pr- and epr-sequences.
He also thanks the anonymous referee, for noting that 
Lemma \ref{same sign lemma} can also be proven using 
Muir's law of extensible minors.


\end{document}